\documentclass[12pt]{article}

\usepackage{amsmath, amsfonts, amssymb}
\usepackage{mathrsfs}
\usepackage{theorem}
\usepackage{bm}

\newtheorem{thm}{Theorem}[section]

\newtheorem{lem}[thm]{Lemma}

{\theorembodyfont{\rmfamily} }
{\theorembodyfont{\rmfamily} }
\newcommand{\norm}[1]{\lVert#1\rVert_{r}^{2}}

\def\R{\mathbb R}

\def\raa{\rangle}

\newcommand{\fin}{\hspace*{\fill}\rule{0.3em}{1ex}}

\numberwithin{equation}{section}

\allowdisplaybreaks

\begin{document}

\title{Well-posedness and Exponential Estimates for the Solutions to Neutral Stochastic Functional Differential Equations  with Infinite Delay
}

\author{Hussein K. Asker ${}^{a,b}$ \footnote{ Email:husseink.askar@uokufa.edu.iq}   \\[0.2cm]
{\small a: Department of Mathematics, Faculty of Computer Science and Mathematics, }\\
{\small Kufa University, Iraq }\\
{\small b:Department of Mathematics, College of science, Swansea University, UK }\\
 }
\maketitle

\begin{abstract}
	
In this work,  neutral stochastic functional differential equations with infinite delay (NSFDEwID) have been addressed. By using the Euler-Maruyama scheme and a localization argument, the existence and uniqueness of solutions to NSFDEwID at the state space $ C_{r} $ under the local weak monotone condition, the weak coercivity condition and the global condition on the neutral term have been investigated. In addition, the $ \mathcal{L}^{2} $ and exponential estimates of NSFDEwID have been studied.

\end{abstract}


\section{Introduction}

 Recently, the neutral stochastic functional differential equations (NSFDEs) have been addressed by many authors, who describe many dynamical systems. Among many other references, we would like to mention \cite{Bao,BYY,hale,hale2,HD,KKMMM,KM,M1,Mo,mw,WC}. Many articles studied wellposedness of NSFDEs by imposing  Lipschitz condition, for example see \cite{bc,Mo,xh,zx}. However, it seemed to be frequently strong in case of the real world. In the last few decades there has been a growing interest in addressing the existence and uniqueness of stochastic functional differential systems under some weaker assumptions. For example \cite{RX,RX2} have studied the existence and uniqueness of solutions with infinite delay at phase space $ BC((- \infty; 0]; \mathbb{R}^{d}) $  with norm $ \lVert \varphi\lVert = \sup_{ - \infty < \theta \le 0} \rvert \varphi (\theta)\rvert $ under non-Lipschitz condition. Bao and Hou \cite{BH}, under a non-Lipschitz condition and a weakened linear growth condition, the existence and uniqueness of mild solutions to stochastic neutral partial functional differential equations have been investigated. Tan et al. \cite{TJS} by the weak convergence approach have reviewed stability in distribution for NSFDEs. von Renesse \cite{von Ren} has used the Euler-Maruyama approximate for SFDEs to show that only weak one-sided local Lipschitz conditions are sufficient for local existence and uniqueness of strong solutions.\\
Based on the approach presented  above, the purpose of this paper is to investigate the existence and uniqueness under the local weak monotone condition, the weak coercivity condition and the global condition on neutral term which can be obtained in the state space $ C_{r} $.\\
The structure of this paper is as follows. Section 2,
recapitulates some basic definitions and notations which has been used to develop our results. Section 3,
gives several sufficient conditions to prove the existence and uniqueness for the equation \eqref{2.2} and the main result. In section 4, the $L^2$
estimate and the exponential estimate have been proved.

\section{Perliminary}

Throughout this paper, unless otherwise specified, we use the following notation. Let  $ R^{d} $ denote the usual $d$-dimensional Euclidean space and   $ \lvert \cdotp \lvert $   the Euclidean norm. If $A$ is a vector or a matrix, its transpose is denoted by $ A^{T} $; and $  \lvert A \rvert = \sqrt{\mbox{trace} (A^{T}A)} $ denotes its trace norm. Denote by $ x^{T}y $ the inner product of  $x$ and $y $ in $\R^{d} $.
Let  $ C((-\infty,0];R^{d}) $ denote the family of all  continuous functions from $ (-\infty,0] $  to $ R^{d}.$
We choose the state space with the fading memory  defined as follows: for given positive number $ r $,
\begin{equation}\label{2.1}
C _{r} = \Big \{ \varphi \in C((-\infty,0];R^{d}) : \|\varphi \| _{r} = \sup_{ - \infty < \theta \leq 0} e^{r\theta} \arrowvert \varphi (\theta) \arrowvert < \infty  \Big \}.
\end{equation}
Then,  $ (C_{r},\|\cdot \| _{r}) $ is a Polish space.  Let $(\Omega,\mathcal{F}, \mathbb{P}) $  be a complete probability space with a filtration $ {\{{\mathcal{F}_{t}\}_{t\in [0,+\infty)}}} $ satisfying the usual conditions (i.e. it is right continuous and $ {{\mathcal{F}_{0}}} $  contains all P-null sets). Let $  1_{B} $ denote the indicator function of a set B. $  M^{2}([ a, b];R^{d})$ is a family of process $ \{ x(t) \}_{a \leq t \leq b}  $ in $ \mathcal{L}^{2}([a, b]; \mathbb{R}^{d}) $ such that $ \mathbb{E} \int_{a}^{b} \rvert x(t)\rvert^{2} dt < \infty $.
Consider a $ d $-dimensional NSFDEwID
\begin{equation}\label{2.2}
d\{x(t) - D(x_{t})\} = b(x_{t})dt + \sigma (x_{t})dw(t), \quad on \quad t \geq 0,
\end{equation}
with the initial data:
\begin{equation}\label{2.3}
x_{0} = \xi = \{\xi(\theta) :-\infty<\theta\leq0 \} \in C_{r},
\end{equation}
where
\begin{equation*}
x_{t} =  x(t+\theta) : - \infty < \theta \leq 0
\end{equation*}
 and $ b , D:  C_{r} \to\mathbb R^{d} $ ;      $ \sigma : C_{r} \to\mathbb R^{d\times m} $ are Borel measurable, $ w(t) $  is an $ m $-dimensional Brownian motion. It should be pointed out that $ x(t) \in  R^{d}  $ is a point and a continuous adapted process $ (x(t))_{t\geq0} $ is called a solution to \eqref{2.2} with the initial value $ x_{0} $, if P-a.s.
  \begin{equation*}
 x (t ) =   \xi (0) + D ( x _{t} ) - D ( \xi) + \int_{0}^{t } b (x_{s}) ds + \int_{0}^{t } \sigma ( x_{s}) dw(s) \quad a.s.,
 \end{equation*}
  while $ x_{t}  \in C_{r}$
 is a continuous function on the interval $ (-\infty , 0] $ taking values in $ R^{d} $ and we call $ (x_{t} )_{t\geq0} $ a functional solution to \eqref{2.2} with the initial value $ x_{0}=\xi \in C_{r}$.

\section{Existence and uniqueness of solutions}

Consider the NSFDEwID \eqref{2.2}, in order to investigate the  existence and uniqueness of  solutions to NSFDEwID, we impose the local weak monotone condition \noindent{\bf (A1)} and the weak coercivity condition \noindent{\bf (A2)}( the condition on the drift coefficient) by developing the tricks adapted in \cite{Mo,M,von Ren,FYM}. Under local weak monotonicity and weak coercivity, the existence and uniqueness for path-independent stochastic differential equations is due to Krylov \cite{Krylov}, which have been extend in [ \cite{r}, Chapter 3]. Under local weak monotone condition and weak coercive condition, \cite{von Ren} studied wellposedness of path-dependent SDEs with finite memory by following Kerylov's approach. The Theorem \ref{mon} extends the result of \cite{von Ren} to NSFDEwID. We assume:

\begin{description}
    \item[(A1)] (local weak monotone condition) There exists a constant  $ L_{R} > 0 $ such that for any $ \phi , \varphi \in C_{r} $ with $ \Arrowvert \varphi \Arrowvert _{r} \vee \Arrowvert \phi \Arrowvert  _{r} \leq R $,
    \begin{equation}\label{A1}
    \begin{split}
    2 \langle \varphi (0) - \phi (0) - \big( D(\varphi) - D (\phi) \big), b ( \varphi) - b ( \phi ) \rangle + \lVert \sigma(\varphi) - \sigma(\phi)\lVert_{r}^{2} \leq L_{R}\lVert \varphi - \phi\lVert_{r}^{2},
    \end{split}
    \end{equation}

    \item[(A2)] (Weak coercivity condition) There exist $ L>0 $ such that:
    \begin{equation}\label{A2}
2\langle \varphi(0)- D(\varphi), b(\varphi) \rangle \vee  |\sigma(\varphi)|^{2}\leq L(1 + \lVert\varphi\lVert_{r}^{2}),
    \end{equation}
    \item[(A3)] There is a $ k \in (0,1) $ such that for all $ \phi , \varphi \in C _{r} $,
    \begin{equation}\label{3.2}
   \rvert D(\phi)-D(\varphi)| \leq k \|\phi - \varphi \Arrowvert_{r} \quad \text{and} \quad D(0)=0.
    \end{equation}
\end{description}
\begin{lem} \label{4.3}
        For all $ \varphi \in C_{r}  $, $ k \in (0, 1) $,
    \begin{equation}
    \rvert \varphi (0) - D(\varphi) \rvert^{2} \leq (1+k)^{2} \lVert \varphi \lVert^{2}_{r}.
    \end{equation}
\end{lem}
\noindent{\bf Proof:} For any $ \varepsilon > 0 $, by the elementary inequality $ \rvert a + b \rvert ^{p} \leq \big [ 1 + \varepsilon ^{\frac{1}{p-1}} \big]^{p-1}  \Big( \rvert a \rvert ^{p} + \dfrac{\rvert b \rvert ^{p}}{\varepsilon} \Big) $, we have:
\begin{equation*}
\rvert \varphi (0) - D(\varphi) \rvert^{2} \leq \big(1+ \varepsilon \big)   \Big( \rvert \varphi (0) \rvert ^{2} + \dfrac{1}{\varepsilon} \rvert D (\varphi) \rvert^{2} \Big).
\end{equation*}
Hence, by \noindent{\bf (A3)}:
\begin{equation*}
\rvert \varphi (0) - D(\varphi) \rvert^{2} \leq \big(1+ \varepsilon  \big)   \Big( \rvert \varphi (0) \rvert ^{2} + \dfrac{k ^{2}}{\varepsilon} \rVert \varphi \rVert^{2} _{r} \Big).
\end{equation*}
Noting that, by letting $ \varepsilon =k $, and
\begin{equation*}
\rvert \varphi (0) \rvert ^{2}  \leq \sup_{ - \infty < \theta \le 0} e ^{2r\theta} \rvert \varphi(\theta) \rvert^{2} = \rVert \varphi \rVert_{r} ^{2},
\end{equation*}
the desired assertion follows.\hfill $\Box$
\begin{thm}\label{mon}
    Under \noindent{\bf (A1)}, \noindent{\bf (A2)} and \noindent{\bf (A3)}, the equation \eqref{2.2} admits a unique strong solution. Moreover, there
    exists constants $ C_2, C_3 > 0 $ such that

  \begin{equation*}
  \mathbb{E}\Big(\sup_{ 0 <  s \le t \wedge \alpha_{R}^{(n)} }e^{2rs}\lVert x^{n}_{s}\lVert ^{2}_{r}\Big)\leq C_3e^{{C_{2}t}},
  \end{equation*}
  where, $ C_{2}= 2C(k\delta)C_{2}(Lk)+C(k\varepsilon_{2}) $ and  
  $ C_{3}=  C(k\delta)\Big[\big(1+ 2(1+k)^{2}\big)\mathbb{E}\norm{\xi}+\dfrac{ 73 L}{r} e^{2rT}\Big]  $.
\end{thm}
\noindent{\bf Proof:} Throughout the whole proof, we assume that $ n,m \geq \frac{r}{\ln2} $ are integers. Define the Euler-Maruyama scheme associated with \eqref{2.2} in the form
\begin{equation}\label{m}
d[x^{n}(t) - D(\hat{x}^{n}_{t})] = b(\hat{x}_{t}^{n})dt + \sigma(\hat{x}^{n}_{t})dw(t), \quad t>0, \quad x^{n}_{0}=\hat{x}^{n}_{0}=x_{0}=\xi,
\end{equation}
where, for each fixed $ t\geq0 $, $ \hat{x}_{t}^{n} \in C_{r} $ is defined in the manner below
\begin{equation}
\hat{x}_{t}^{n}(\theta): = x^{n}((t+\theta) \wedge t_{n}), \quad \theta \in (- \infty,0], \quad t_{n}:= \frac{\left \lfloor{nt}\right \rfloor}{n}.
\end{equation}
We point out that all constants are independent of $ n\geq 1 $. By the fixed point theorem, the path-dependent NSFDE \eqref{m} has a unique solution by solving pice-wisely with the time step lenght $ \frac{1}{n} $. For any $ R > 3 \norm{\xi} $, define the stopping times
\begin{equation}\label{stop}
\tau_{R}^{(n)}:=\inf \big\{ t \geq 0:\rvert x^{n}(t)\rvert \geq \frac{R}{3}\} \quad \and \quad \alpha_{R}^{(n)}:=\inf \big\{ t \geq 0:\lVert x^{n}_{t}\lVert _{r}\geq \frac{R}{3}\}.
\end{equation}
by [Wu, \cite{FYM}], we have $ \tau_{R}^{(n)}=\alpha_{R}^{(n)} $. Observe that
\begin{equation*}
\lVert x^{n}_{t}\lVert _{r}= \sup_{ - \infty < \theta \le 0}(e^{r\theta} \rvert x^{n}(t+\theta)\rvert ) \geq e ^{r(t_{n}- t)} \rvert x^{n}(t_{n})\rvert,
\end{equation*}
which, in addition to $ n \geq \frac{r}{\ln2} $, we have $ e^{\frac{r}{n}}\leq2 $ yields that
\begin{equation*}
\rvert x^{n}(t_{n})\rvert \leq e^{r(t-t_{n})} \lVert x^{n}_{t}\lVert_{r}= e^{r(t-\frac{\left \lfloor{nt}\right \rfloor}{n})} \lVert x^{n}_{t}\lVert_{r}= e^{r(\frac{nt-\left \lfloor{nt}\right \rfloor}{n})} \lVert x^{n}_{t}\lVert_{r} \leq e^{\frac{r}{n}}\lVert x^{n}_{t}\lVert _{r} \leq 2\lVert x^{n}_{t}\lVert _{r}.
\end{equation*}
This, together with
\begin{equation*}
\begin{split}
\lVert \hat{x}^{n}_{t}\lVert _{r}&=\sup_{ - \infty < \theta \le 0}(e^{r\theta}\hat{x}^{n}_{t}(\theta)) =\sup_{ - \infty < \theta \le 0}(e^{r\theta}x^{n}((t+\theta)\wedge t_{n})\\
& = \sup_{ - \infty < \theta \le 0}(e^{r\theta}x^{n}((t+\theta)\wedge t_{n})1_{\{t+\theta \leq t_{n}\}} +\sup_{ - \infty < \theta \le 0}(e^{r\theta}x^{n}((t+\theta)\wedge t_{n})1_{\{t+\theta \geq t_{n}\}}\\
&  = \sup_{ - \infty < \theta \le 0}(e^{r\theta}x^{n}(t+\theta)) +\sup_{ - \infty < \theta \le 0}(e^{r\theta}x^{n}(t_{n}))\\
& \leq \lVert x^{n}_{t}\lVert _{r}+ \rvert x^{n}(t_{n})\rvert,
\end{split}
\end{equation*}
further leads to
\begin{equation}\label{Loc1}
\lVert \hat{x}^{n}_{t}\lVert _{r} \leq 3\lVert x^{n}_{t}\lVert _{r}.
\end{equation}
Let
\begin{equation*}
C(R):= \sup_{ \lVert\zeta\lVert _{r}\leq \frac{R}{3}}\rvert b(\zeta)\rvert < \infty.\quad \text{So},\quad \lvert b(x_{t}^{n})\lvert \leq C(R)  \quad R\in(3\norm{\xi}, \infty),\quad t\in[0, \tau_{R}^{(n)}].
\end{equation*}
Then, it follows from \eqref{Loc1} and the notation of $ \tau_{R}^{(n)} $ that
\begin{equation}\label{Loc2}
\rvert b(\hat{x}_{t}^{n})\rvert\leq C(R), \quad t\leq \tau_{R}^{(n)}=\alpha_{R}^{(n)}.
\end{equation}
To show that $ (x^{n}_{t})_{t\geq0} $ converges in probability to some stochastic process $ (x_{t})_{t\geq0} $ as $ n \rightarrow \infty$, for $ n,m \geq \frac{r}{\ln2} $, 
 set $ z^{n,m}(t) := x^{n}(t) - x^{m}(t) $, $ D^{n,m}(\hat{x}_{t}):= D(\hat{x}_{t}^{n}) -D(\hat{x}_{t}^{m}) $, $ \hat{z}^{n,m}(t) := \hat{x}^{n}(t) - \hat{x}^{m}(t) $,  $ p^{n}_{t}:= x^{n}_{t} - \hat{x}^{n}_{t} $ and $ \varGamma^{n,m}(t):= z^{n,m}(t) - D^{n,m}(\hat{x}_{t})  $, by using the fact that $ x^{n}_{t} $ and $  x^{m}_{t} $ share the initial value with the elementary inequality and \noindent{\bf (A3)},  we note that for any $ \varepsilon_{1} = \frac{k}{1-k}   $ and  $ \varepsilon_{2} > \frac{k}{1-k}   $, we have
\begin{equation*}
\begin{split}
e^{2rt} \norm{z^{n,m}_{t}} & = \sup_{ 0 <  s \le t} (e^{2rs}\rvert z^{n,m}(s)\rvert^{2}) \\
&= \sup_{ 0 <  s \le t} (e^{2rs}\rvert z^{n,m}(s) - D^{n,m}(\hat{x}_{s}) + D^{n,m}(\hat{x}_{s})\rvert^{2})\\
& \leq \dfrac{1}{1-k} \sup_{ 0 <  s \le t} (e^{2rs}\rvert z^{n,m}(s) - D^{n,m}(\hat{x}_{s})\rvert^{2}) +k e^{2rt} \norm{\hat{x}^{n}_{t}-\hat{x}^{m}_{t}}, \\
& \leq \dfrac{1}{1-k} \sup_{ 0 <  s \le t} (e^{2rs}\rvert \varGamma^{n,m}(s) \rvert^{2}) +k(1+\varepsilon_{2}) e^{2rt} \norm{p^{n}_{t}-p^{m}_{t}}  \\
&\qquad+\dfrac{k(1+\varepsilon_{2})}{\varepsilon_{2}} e^{2rt} \norm{x^{n}_{t}-x^{m}_{t}},\\
\end{split}
\end{equation*}
note that $ \delta:=\dfrac{k(1+\varepsilon_{2})}{\varepsilon_{2}} < 1  $, therefore
\begin{equation}\label{Loc3}
e^{2rt} \norm{z^{n,m}_{t}}  \leq \dfrac{ 1}{(1- k)(1-\delta)}\sup_{ 0 <  s \le t} (e^{2rs}\rvert \varGamma^{n,m}(s)\rvert^{2})+\dfrac{k(1+\varepsilon_{2})}{1-\delta} e^{2rt} \norm{p^{n}_{t}-p^{m}_{t}} .\\
\end{equation}
Also, with the definition of $ \tau_{R}^{(n)} $ and $ \alpha_{R}^{(n)} $, \eqref{Loc1} implies that
\begin{equation}\label{p1}
\rVert p_{t}^{n}\rVert_{r}\leq \frac{4}{3}R, \quad t\leq \tau_{R}^{(n)}
\end{equation}
Now, for $ t\leq \tau_{R}^{(n)}\wedge \tau_{R}^{(m)}=\alpha_{R}^{(n)}\wedge \alpha_{R}^{(m)} $, by the It\^o formula, we derive from \noindent{\bf (A1)}, \eqref{Loc2},\eqref{Loc3} \eqref{p1} and \noindent{\bf (A3)}   that
\begin{equation}
\begin{split}
&d(e^{2rt} \rvert \varGamma^{n,m}(t) \rvert^{2})=e^{2rt}\{ 2r \rvert \varGamma^{n,m}(t)\rvert^{2} + 2 \langle \varGamma^{n,m}(t), b(\hat{x}^{n}_{t}) - b(\hat{x}^{m}_{t}) \rangle\\
& \quad + \lVert\sigma(\hat{x}^{n}_{t})- \sigma (\hat{x}^{m}_{t}) \lVert^{2}_{r} \}dt + dM^{n,m}(t) \\
& = e^{2rt}\{ 2r \rvert \varGamma^{n,m}(t)\rvert^{2} + 2 \langle z^{n,m}(t) - D^{n,m}(\hat{x}_{t}) + \hat{z}^{n,m}(t)- \hat{z}^{n,m}(t) , b(\hat{x}^{n}_{t}) - b(\hat{x}^{m}_{t}) \rangle\\
& \quad+ \lVert\sigma(\hat{x}^{n}_{t})- \sigma (\hat{x}^{m}_{t}) \lVert^{2}_{r} \}dt + dM^{n,m}(t)\\
& = e^{2rt}\{ 2r \rvert \varGamma^{n,m}(t)\rvert^{2} + 2 \langle z^{n,m}(t)- \hat{z}^{n,m}(t),b(\hat{x}^{n}_{t}) - b(\hat{x}^{m}_{t}) \rangle + 2 \langle \hat{z}^{n,m}(t) - D^{n,m}(\hat{x}_{t}),\\
&\qquad b(\hat{x}^{n}_{t}) - b(\hat{x}^{m}_{t}) \rangle+ \lVert\sigma(\hat{x}^{n}_{t})- \sigma (\hat{x}^{m}_{t}) \lVert^{2}_{r} \}dt +dM^{n,m}(t)\\
&  \leq e^{2rt}\{ 2r \rvert \varGamma^{n,m}(t)\rvert^{2}  + 2\langle  p^{n}_{t}(0) - p^{m}_{t}(0)  ,  b(\hat{x}^{n}_{t}) - b(\hat{x}^{m}_{t})\rangle + L_{R}\norm{\hat{x}^{n}_{t}- \hat{x}^{m}_{t}}\} dt\\
& \qquad +dM^{n,m}(t)\\
&  \leq e^{2rt}\{ 2r \rvert \varGamma^{n,m}(t)\rvert^{2}  + 2\rvert  p^{n}_{t}(0) - p^{m}_{t}(0) \rvert \cdot\big( \rvert b(\hat{x}^{n}_{t})\rvert +\rvert b(\hat{x}^{m}_{t})\rvert\big)\\
&\qquad +2L_{R}\norm{p^{n}_{t}-p^{m}_{t}} + 2L_{R}\norm{z^{n,m}_{t}}\} dt+dM^{n,m}(t)\\
&  \leq e^{2rt}\Big\{ \Big(2r+\dfrac{ 2L_{R}}{(1- k)(1-\delta)}\Big) \rvert \varGamma^{n,m}(t)\rvert^{2} + 4C(R) \rvert  p^{n}_{t}(0)\rvert + 4C(R) \rvert p^{m}_{t}(0)\rvert \\
& \qquad+2L_{R}(1+\dfrac{k(1+\varepsilon_{2})}{1-\delta} )\norm{p^{n}_{t}-p^{m}_{t}}\Big\} dt+dM^{n,m}(t)\\
&  \leq  K \{ \sup_{ 0 <  s \le t} (e^{2rs} \rvert \varGamma^{n,m}(s)\rvert^{2} + e^{2rt} \big(\lVert  p^{n}_{t}\lVert_{r} + \lVert p^{m}_{t}\lVert_{r}\big)   \} dt + dM^{n,m}(t) ,
\end{split}
\end{equation}
where $ K\geq \Big(2r+\dfrac{2L_{R}}{(1- k)(1-\delta)}\Big) + \Big(4C(R)+4L_{R}\big(1+\dfrac{k(1+\varepsilon_{2})}{1-\delta} \big)\Big)>0$, $ dM^{n,m}(t)=2e^{2rt} \langle \varGamma^{n,m}(t), (\sigma(\hat{x}^{n}_{t})- \sigma (\hat{x}^{m}_{t})) dw(t) \raa $. Thus, for fixed $ T >0 $, each $ p \in (0, 1) $ and $ \alpha > \frac{1+p}{1-p} $, in terms of stochastic Gronwall inequality [\cite{von Ren}, Lemma 5.4], there exist constant $ c>0 $ depending on $ p, \alpha $ such that

\begin{equation}\label{3.8}
\begin{split}
\mathbb{E} & \Big(  \sup_{ 0 \leq t \leq T \wedge  \tau_{R}^{(n)}\wedge \tau_{R}^{(m)}} e^{2rt} \rvert \varGamma^{n,m}(s)\rvert^{2}\Big)^{p}  \leq \\
& \qquad\quad c\Big(  \int_{0}^{T} \mathbb{E} (\lVert p^{n}_{t}\lVert_{r}^{\alpha}1_{\{ t \leq \tau_{R}^{(n)} \}} )dt\Big)^{\frac{p}{\alpha}} + c\Big (\int_{0}^{T} \mathbb{E} (\lVert p^{m}_{t}\lVert_{r}^{\alpha}1_{\{ t \leq \tau_{R}^{(m)} \}} )dt\Big)^{\frac{p}{\alpha}}.
\end{split}
\end{equation}
Finally, to estimate $ p^{n}_{ t}  $, for $ \theta \leq 0 $, note that \eqref{m} implies
\begin{equation*}
\begin{split}
p^{n}_{t}(\theta) & = x^{n}(t+\theta) - \hat{x}^{n}(t+\theta) =x^{n}(t+\theta) -x^{n}((t+\theta) \wedge t_{n})
\\
&  = \left\{ \begin{array}{rcl}
0, \qquad \qquad \qquad \qquad & \mbox{for}
& t+\theta \leq t_{n} ,\\
D(\hat{x}^{n}_{t+\theta})-D(\hat{x}^{n}_{t_{n}})+ \int_{t_{n}}^{t+\theta}b(\hat{x}^{n}_{s})ds + \int_{t_{n}}^{t+\theta}\sigma(\hat{x}^{n}_{s})dw(s),  & \mbox{for} & t+\theta > t_{n},
\end{array}\right.
\end{split}
\end{equation*}
in view of \noindent{\bf (A3)}, we deduce that
\begin{equation}
\begin{split}
&\lVert p^{n}_{t}\lVert_{r} = \sup_{ - \infty < \theta \le 0}(e^{r\theta}\rvert p^{n}_{t}(\theta)\rvert)= \sup_{ - \infty < \theta \le 0}(e^{r\theta}\rvert x^{n}(t+\theta) -x^{n}((t+\theta)\wedge t_{n})\rvert)1_{\{t+\theta\geq t_{n}\}}\\
&\qquad = e^{-rt}\sup_{ t_{n} < s \le t}(e^{rs}\rvert x^{n}(s) -x^{n}(s\wedge t_{n})\rvert)\\
& \leq e^{-rt} \sup_{ t_{n} \leq s \le t} \Big(e^{rs}\rvert D(\hat{x}^{n}_{s})-D(\hat{x}^{n}_{s\wedge t_{n}})\rvert+ e^{rs} \int_{t_{n}}^{s}\rvert b(\hat{x}^{n}_{u})\rvert du + e^{rs}\Big \rvert \int_{t_{n}}^{s}\sigma(\hat{x}^{n}_{u})dw(u)\Big \rvert  \Big)\\
& \leq ke^{-rt} \sup_{ t_{n} \leq s \le t}\Big(e^{rs}\lVert\hat{x}^{n}_{s}-\hat{x}^{n}_{ s\wedge t_{n}} \lVert_{r}\Big)+ \int_{t_{n}}^{t}\rvert b(\hat{x}^{n}_{s})\rvert ds  + \sup_{ t_{n} \leq s \le t, }\Big \rvert \int_{t_{n}}^{s}\sigma(\hat{x}^{n}_{u})dw(u)\Big \rvert.
\end{split}
\end{equation}
Note that, 
\begin{equation*}
\begin{split}
\lVert\hat{x}^{n}_{s}-\hat{x}^{n}_{ s\wedge t_{n}} \lVert_{r}&= \sup_{ -\infty \leq \theta \le 0}\Big(e^{r\theta}\rvert\hat{x}^{n}(s+\theta)-\hat{x}^{n} ((s+\theta)\wedge t_{n}) \rvert\Big)1_{\{t_{n}\leq s+\theta\}}\\
& \leq e^{-rs} \sup_{ t_{n} \leq u \le s}\big(e^{ru}\rvert\hat{x}^{n}(u)-\hat{x}^{n} (t_{n}) \rvert\Big)=0,
\end{split}
\end{equation*}
so that,
\begin{equation}
\begin{split}\label{p}
&\lVert p^{n}_{t}\lVert_{r}
\leq  \int_{t_{n}}^{t}\rvert b(\hat{x}^{n}_{s})\rvert ds  + \sup_{ t_{n} \leq s \le t, }\Big \rvert \int_{t_{n}}^{s}\sigma(\hat{x}^{n}_{u})dw(u)\Big \rvert.\\
\end{split}
\end{equation}
By means of \eqref{Loc2}, for any $ t\geq0 $ and $ t\leq \tau_{R}^{(n)}=\alpha_{R}^{(n)} $, we infer that
\begin{equation}\label{b}
\lim\limits_{n\rightarrow\infty}\mathbb{E}\int_{t_{n}}^{t\wedge \tau_{R}^{(n)}}\rvert b(\hat{x}^{n}_{s})\rvert ds\leq\lim\limits_{n\rightarrow\infty} \dfrac{1}{n} C(R) = 0 .
\end{equation}
On the other hand,  the Burkhold-Davis-Gundy inequality, together with the local boundedness of $ \sigma $,  there exists a constant $ M_{R} > 0 $ such that, for any $ t\geq 0  $ and $ t\leq \tau_{R}^{(n)}=\alpha_{R}^{(n)} $,
\begin{equation*}
\lim\limits_{n\rightarrow\infty}\mathbb{E} \Big(\sup_{t_{n} \leq s \le t \wedge \tau_{R}^{(n)}} \Big\rvert \int_{t_{n}}^{s}\sigma(\hat{x}^{n}_{s})dw(s)\Big\rvert^{2} \Big)\leq   \lim\limits_{n\rightarrow\infty} \dfrac{M_{R}}{n}=0,
\end{equation*}
thus, from this with \eqref{p} and \eqref{b}, 
we conclude that
\begin{equation}\label{3.9}
\sup_{ t \in [0, T]}\mathbb{E}(\lVert p^{n}_{t}\lVert_{r}^{\alpha}1_{\{ t \leq \tau_{R}^{(n)}\}}) \rightarrow 0 \quad as \quad n \rightarrow \infty.
\end{equation}
Subsequently, for $ p=\dfrac{1}{2} $ taking \eqref{3.8} and \eqref{Loc3} into consideration implies that
\begin{equation}\label{3.10}
\lim\limits_{n,m \rightarrow \infty}\mathbb{P}\Big\{\sup_{ 0 \leq t \leq T \wedge  \tau_{R}^{(n)}\wedge \tau_{R}^{(m)}} \lVert x^{n}_{t} - x^{m}_{t} \lVert_{r} \}=0.
\end{equation}
Next, to ensure that $ x^{n} $ · converges in probability to a solution of \eqref{2.2}, it remains to prove
\begin{equation}\label{3.11}
\lim\limits_{R \rightarrow \infty} \limsup_{n \rightarrow \infty}\mathbb{P}( \tau_{R}^{(n)} \leq T) = 0.
\end{equation}
Indeed, \eqref{3.10} and \eqref{3.11} yield that
\begin{equation*}
\lim\limits_{n,m \rightarrow \infty}\mathbb{P}\Big\{\sup_{ 0 \leq t \leq T } \lVert x^{n}_{t} - x^{m}_{t} \lVert_{r} \geq \varepsilon\}=0, \quad \varepsilon > 0,
\end{equation*}
Whence, by keeping in mind that $ C_{r} $ is a Polish space under the metric $ \lVert \cdot\lVert_{r}$( completeness of $ (C_{r},\lVert \cdot\lVert_{r}) $), there exists a continuous adapted process stochastic process $ (x_{t})_{t\in[0, T]} $ on $ C_{r} $ such that
\begin{equation}\label{x}
\sup_{ 0 \leq t \leq T } \lVert x^{n}_{t} - x_{t} \lVert_{r} \rightarrow 0 \quad \text{in} \quad \text{probability}\quad \text{as}  \quad n \rightarrow \infty.
\end{equation}
Then, by following the standard argument we can show that $ (x_{t})_{t\in[0, T]} $ is the unique functional solution to \eqref{2.2} under \noindent{\bf (A1)}, \noindent{\bf (A2)} and \noindent{\bf (A3)}.

So, in what follows, it remains to show that \eqref{3.11} holds true. Set $ \Gamma^{n}(t)=x^{n}(t)-D(\hat{x}^{n}_{t}) $.
Note, by the elementary inequality and \noindent{\bf (A3)} similar to \eqref{Loc3} with $ x^{n}_{0}=\xi $, for $ \varepsilon > \frac{3k^{2}}{1- 3k^{2}}$,  that
\begin{equation*}
\begin{split}
e^{2rt}\norm{x^{n}_{t}}&\leq \norm{\xi}+(1+\varepsilon) \sup_{ 0 <  s \le t}  \Big(e^{2rs}\rvert \Gamma^{n}(s)\rvert^{2}\Big)+\dfrac{3k^{2}(1+\varepsilon)}{\varepsilon} e^{2rt} \norm{x^{n}_{t}},
\end{split}
\end{equation*}
set, $\gamma:=\dfrac{3k^{2}(1+\varepsilon)}{\varepsilon} < 1 $, thus
\begin{equation}\label{n}
\begin{split}
e^{2rt}\norm{x^{n}_{t}}&\leq\dfrac{1}{1-\gamma} \norm{\xi}+\dfrac{(1+\varepsilon)}{1-\gamma} \sup_{ 0 <  s \le t}  \Big(e^{2rs}\rvert \Gamma^{n}(s)\rvert^{2}\Big).
\end{split}
\end{equation}
By It\^o's formula, we deduce from  Lemma \ref{4.3}, \noindent{\bf (A3)},  \eqref{Loc1}, \noindent{\bf (A2)}, and \eqref{Loc2}  that
\begin{equation}\label{3.14}
\begin{split}
&e^{2rt}\rvert\Gamma^{n}(t)\rvert^{2} =  \rvert \xi(0) - D ( \xi) \rvert^{2} +   \int_{0}^{t} e ^{2rs} \Big [ 2r \rvert \Gamma^{n}(s) \rvert^{2} \\
&  \qquad + 2\langle x^{n}(s) - D (\hat{x}^{n}_{s}), b(\hat{x}^{n}_{s })\rangle + \lVert \sigma (\hat{x}^{n}_{s}) \lVert_{r}^{2} \Big] ds + dN^{n}(t)\\
& \leq   (1+k)^{2}\norm{\xi} +   \int_{0}^{t} e ^{2rs} \Big [ 4r \rvert x^{n}(s)\rvert^{2} + 4r k^{2} \norm{\hat{x}^{n}_{s}}+2\langle x^{n}(s) - D (\hat{x}^{n}_{s}) + \hat{x}^{n}(s)- \hat{x}^{n}(s)\\
&\qquad b(\hat{x}^{n}_{s })\rangle  +  \lVert \sigma (\hat{x}^{n}_{s}) \lVert_{r}^{2} \Big] ds + dN^{n}(t)\\
& \leq   (1+k)^{2}\norm{\xi} +   \int_{0}^{t} e ^{2rs} \Big [ 4r \norm{x^{n}_{s}} + 12rk^{2}\norm{x^{n}_{s}}  + 2\langle x^{n}(s) - \hat{x}^{n}(s), b(\hat{x}^{n}_{s }) \raa   \\
&  \qquad+ 2\langle \hat{x}^{n}(s)- D (\hat{x}^{n}_{s}), b(\hat{x}^{n}_{s })\raa +  \lVert \sigma (\hat{x}^{n}_{s}) \lVert_{r}^{2} \Big] ds + dN^{n}(t)\\
& \leq   (1+k)^{2}\norm{\xi} +   \int_{0}^{t} e ^{2rs} \Big [4r (1+ 3k^{2}) \norm{x^{n}_{s}}  + 2 \rvert p^{n}_{s}(0)\rvert\rvert b(\hat{x}^{n}_{s })\rvert \\
& \qquad+  L(1+ \norm{\hat{x}^{n}_{s}}) \Big] ds + dN^{n}(t)\\
& =   (1+k)^{2}\norm{\xi}+ \dfrac{L}{2r}e ^{2rt} + C_{1}(Lk)  \int_{0}^{t} e ^{2rs}   \norm{x^{n}_{s}} ds  + 2C(R) \int_{0}^{t}e ^{2rs}\rvert p^{n}_{s}(0)\rvert ds+dN^{n}(t) ,
\end{split}
\end{equation}
where, $ C_{1}(Lk)=\big(4r (1+ 3k^{2})+3 L \big) $ and $ dN^{n}(t)=2\int_{0}^{t} e ^{2rs}  \langle   \Gamma^{n}(s) ,\sigma (\hat{x}^{n}_{s}) dw(s)\rangle $.
  Thus for any $ t \in [0, T], $ from \eqref{3.14} we have,
\begin{equation}\label{3.15}
\begin{split}
\mathbb{E} \Big(\sup_{ 0 <  s \le t \wedge  \tau_{R}^{(n)}}&e^{2rs}\rvert \Gamma^{n}(s)\rvert^{2}\Big) \leq (1+k)^{2}\mathbb{E}\norm{\xi}+\dfrac{  L}{2r} e^{2rT} \\
& +C_{1}(Lk)\mathbb{E}\int_{0}^{t}e ^{ 2rs} \lVert x^{n}_{s}\lVert ^{2}_{r}ds + 2 C(R)\int_{0}^{t}e ^{2rs} \mathbb{E} \Big(\rvert p _{s}^{n}(0)\rvert 1_{\{ s \leq \tau_{R}^{(n)}\}}\Big) ds\\
&\qquad \qquad + 2 \mathbb{E}\Big(\sup_{ 0 \leq s \leq t \wedge  \tau_{R}^{(n)}}\int_{0}^{s} e ^{2ru}  \langle   \Gamma^{n}(u) ,\sigma (\hat{x}^{n}_{u})dw(u)\rangle\Big).\\
\end{split}
\end{equation}
Next, by the Burkholder - Davis - Gundy inequality, we infer that
\begin{equation}\label{3.16b}
 \begin{split}
 2 \mathbb{E}&\Big(\sup_{ 0 \leq s \leq t \wedge  \tau_{R}^{(n)}}\int_{0}^{s} e ^{2ru}  [   \Gamma^{n}(u) ]^{T}\sigma (\hat{x}^{n}_{u})dw(u)\Big)\\
 & \leq 8\sqrt{2} \mathbb{E} \Big ( \int_{0}^{t \wedge  \tau_{R}^{(n)}} e ^{2 ru} \rvert [\Gamma^{n}(u)] ^{T} \sigma (\hat{x}^{n}_{u})\rvert^{2} du \Big) ^{\frac{1}{2}}\\
 & \leq 12 \mathbb{E} \Big ( \int_{0}^{t \wedge  \tau_{R}^{(n)}} e ^{2 ru} \rvert \Gamma^{n}(u) \rvert ^{2}  \lvert \sigma (\hat{x}^{n}_{u})\lvert^{2} du \Big)^{\frac{1}{2}}\\
 & \leq 12 \mathbb{E} \Big[ \Big ( \sup_{ 0 <  u \le t \wedge  \tau_{R}^{(n)}}  e ^{ 2ru} \rvert \Gamma^{n}(u) \rvert ^{2} \Big)^{\frac{1}{2}} \Big( \int_{0}^{t \wedge  \tau_{R}^{(n)}} e ^{ 2ru} \lvert \sigma (\hat{x}^{n}_{u})\lvert^{2} du \Big)^{\frac{1}{2}} \Big]\\
 & \leq \dfrac{1}{2} \mathbb{E} \Big ( \sup_{ 0 \leq  u \le t \wedge \tau_{R}^{(n)} }  e ^{ 2ru} \rvert \Gamma^{n}(u) \rvert ^{2} \Big)+ 72 \mathbb{E}  \int_{0}^{t \wedge\tau_{R}^{(n)}} e ^{ 2ru}  \lvert \sigma (\hat{x}^{n}_{u})\lvert^{2} du\\
 &\leq \dfrac{1}{2} \mathbb{E} \Big ( \sup_{ 0 <  u \le t \wedge \tau_{R}^{(n)} }  e ^{ 2ru} \rvert \Gamma^{n}(u) \rvert ^{2} \Big) + 72L \mathbb{E}  \int_{0}^{t \wedge\tau_{R}^{(n)}} e ^{ 2ru}   (1+ \norm{\hat{x}^{n}_{u}}) du\\
 & \leq \dfrac{1}{2} \mathbb{E} \Big ( \sup_{ 0 <  u \le t \wedge \tau_{R}^{(n)} }  e ^{ 2ru} \rvert \Gamma^{n}(u) \rvert ^{2} \Big) + 72 L\mathbb{E}  \int_{0}^{t \wedge\tau_{R}^{(n)}} e ^{ 2ru}   (1+ 3\norm{x^{n}_{u}}) du\\
& \leq \dfrac{1}{2} \mathbb{E} \Big ( \sup_{ 0 <  u \le t \wedge \tau_{R}^{(n)} }  e ^{ 2ru} \rvert \Gamma^{n}(u) \rvert ^{2} \Big) +\dfrac{36L}{r}e ^{ 2rT}  + 216L\mathbb{E}  \int_{0}^{t }e ^{ 2ru} \lVert x^{n}_{u}\lVert ^{2}_{r}du.\\
 \end{split}
\end{equation}
Plugging \eqref{3.16b} into \eqref{3.15}, we get that
\begin{equation}
\begin{split}
\mathbb{E} \Big(\sup_{ 0 <  s \le t \wedge  \tau_{R}^{(n)}}&e^{2rs}\rvert \Gamma^{n}(s)\rvert^{2}\Big) \leq (1+k)^{2}\mathbb{E}\norm{\xi}+\dfrac{ 73 L}{2r} e^{2rT} \\
& +C_{2}(Lk)\mathbb{E}\int_{0}^{t}e ^{ 2rs} \lVert x^{n}_{s}\lVert ^{2}_{r}ds + 2 C(R)\int_{0}^{t}e ^{2rs} \mathbb{E} \Big(\rvert p _{s}^{n}(0)\rvert 1_{\{ s \leq \tau_{R}^{(n)}\}}\Big) ds\\
&\qquad \qquad \dfrac{1}{2} \mathbb{E} \Big ( \sup_{ 0 <  s \le t \wedge \tau_{R}^{(n)} }  e ^{ 2rs} \rvert \Gamma^{n}(s) \rvert ^{2} \Big) ,\\
\end{split}
\end{equation}
where, $ C_{2}(Lk)=C_{1}(Lk)+216L $. Consequently,
\begin{equation}\label{3.19}
\begin{split}
\mathbb{E} \Big(\sup_{ 0 <  s \le t \wedge  \tau_{R}^{(n)}}e^{2rs} \rvert \Gamma^{n}(s)\rvert^{2}\Big)& \leq 2 (1+k)^{2}\mathbb{E}\norm{\xi}+\dfrac{ 73 L}{r} e^{2rT} \\
& +2C_{2}(Lk)\mathbb{E}\int_{0}^{t}e ^{ 2rs} \lVert x^{n}_{s}\lVert ^{2}_{r}ds \\
& + 4 C(R)\int_{0}^{t}e ^{2rs} \mathbb{E} \Big(\rvert p _{s}^{n}(0)\rvert 1_{\{ s \leq \tau_{R}^{(n)}\}}\Big) ds .
\end{split}
\end{equation}
Considering \eqref{n}, the \eqref{3.19} gives that, for any $ t \in [0, T] $, 
\begin{equation}
\begin{split}
&\mathbb{E} \Big(\sup_{ 0 <  s \le t \wedge  \tau_{R}^{(n)}}e^{2rs} \rvert \Gamma^{n}(s)\rvert^{2}\Big)\leq \Big[\Big(2(1+k)^{2}+ \dfrac{2C_{2}(Lk)t}{1-\gamma}\Big)\mathbb{E}\norm{\xi}+\dfrac{ 73 L}{r} e^{2rT}\\
& +4 C(R)\int_{0}^{t}e ^{2rs} \mathbb{E} \Big(\rvert p _{s}^{n}(0)\rvert 1_{\{ s \leq \tau_{R}^{(n)}\}}\Big) ds \Big] + C_{2}\mathbb{E}\int_{0}^{t} \Gamma^{n}(s)ds,
\end{split}
\end{equation}
 So, the Gronwall inequality yields that
\begin{equation}\label{3.17}
\Gamma^{n, R}(t) \leq C_{1}e^{C_{2}T},
\end{equation}
where, $ C_{2}= \dfrac{2C_{2}(Lk)(1+\varepsilon) }{1-\gamma}$, and
 \begin{equation*}
  \begin{split}
  C_{1}=  \Big[\Big(2(1+k)^{2}+ \dfrac{2C_{2}(Lk)t}{1-\gamma}\Big)\mathbb{E}\norm{\xi}+\dfrac{ 73 L}{r} e^{2rT}+4 C(R)\int_{0}^{t}e ^{2rs} \mathbb{E} \Big(\rvert p _{s}^{n}(0)\rvert 1_{\{ s \leq \tau_{R}^{(n)}\}}\Big) ds \Big]
  \end{split}
 \end{equation*} 
According to the notion of $ \tau_{R}^{(n)} $, the event
\begin{equation*}
\Big\{\tau_{R}^{(n)}\leq T, \sup_{ 0 \leq t \le T \wedge \tau_{R}^{(n)}} \rvert x^{n}(t)\rvert < \dfrac{R}{4}\Big\}
\end{equation*}
is empty set, which, together with \eqref{3.9}, \eqref{3.17} and Chebyshev's inequality yields that
\begin{equation}\label{3.18}
\begin{split}
\lim\limits_{R\rightarrow\infty}\lim\limits_{n\rightarrow\infty}\mathbb{P}(\tau_{R}^{(n)} \leq T) &=\lim\limits_{R\rightarrow\infty}\lim\limits_{n\rightarrow\infty} \mathbb{P}\Bigg(\tau_{R}^{(n)} \leq T,\sup_{ 0 \leq t \le \tau_{R}^{(n)}\wedge T} \rvert x^{n}(t)\rvert \geq \dfrac{R}{4} \Bigg)\\
&\leq\lim\limits_{R\rightarrow\infty}\lim\limits_{n\rightarrow\infty} \mathbb{P}\Bigg(\sup_{ 0 \leq t \le \tau_{R}^{(n)}\wedge T} \rvert x^{n}(t)\rvert \geq \dfrac{R}{4} \Bigg)\\
 & \leq\lim\limits_{R\rightarrow\infty}\lim\limits_{n\rightarrow\infty} \dfrac{16\Gamma^{n, R}(T) }{R^{2}}=0.
\end{split}
\end{equation}
So, \eqref{3.11} holds.
Finally, by \eqref{3.9}, \eqref{3.10} and \eqref{3.17} and employing Fatou’s lemma for $ n \rightarrow \infty $, we obtain there exists constants $ C_3, C_4 > 0 $ such that
\begin{equation*}
\mathbb{E}\Big(\sup_{ 0 <  s \le t \wedge \alpha_{R}^{(n)} }e^{2rs}\lVert x^{n}_{s}\lVert ^{2}_{r}\Big)\leq C_3e^{{C_{2}t}},
\end{equation*}
where $ C_{3}=  C(k\delta)\Big[\big(1+ 2(1+k)^{2}\big)\mathbb{E}\norm{\xi}+\dfrac{ 73 L}{r} e^{2rT}\Big]  $
 and the proof complete. \hfill $\Box$\\
\section{Exponential estimate} 
Let $ x(t)$ where $t\in [0,\infty) $, be a unique solution to NSFDEwID \eqref{2.2}. In this section,  first we derive the $L^2$
estimate, then we obtain the exponential estimate for the solution.
\begin{thm}\label{est}
	Under  \noindent{\bf (A2)} and \noindent{\bf (A3)}, for  $ t\geq0 $ there exists constants 	 $ C_{4} ,C_5> 0 $ such that the solution $ x(t) $ of \eqref{2.2}, satisfies: 
	\begin{equation}
	E[\sup_{0\leq s\leq t}|x(s)|^2]\leq C_{4}e^{C_{5}t},
	\end{equation}
	where $ C_{4}= \dfrac{1}{1-2k^{2}}\Big[\Big( 2k^{2}e ^{- 2r t}  + 4(1+k)^{2}+\dfrac{146L}{r}\Big) \mathbb{E} \rVert  \xi \rVert^{2}_{r}+ 292 LT\Big] $ and $ C_{5}=292 L $.
\end{thm}
\textbf{Proof:} By using the inequality $ (a+b)^{2}\leq 2a^{2}+2b^{2}$ and $ $\noindent{\bf (A3)}, we have
\begin{equation}\label{ee}
\begin{split}
\mathbb{E}\big(\sup_{0\leq s\leq
	t} \rvert x(s) \rvert^{2}\big)&\leq 2\mathbb{E}  \Big ( \sup_{ 0 < s\le t}
\rvert D ( x_{s})\rvert^{2} \Big)+2\mathbb{E}  \Big ( \sup_{ 0 < s\le t}
\rvert x(s) - D ( x_{s})\rvert^{2} \Big) \\
& \leq 2k^{2}E\|x_{t}\|_r^2+ 2\mathbb{E}  \Big ( \sup_{ 0 < s\le t}
\rvert x(s) - D ( x_{s})\rvert^{2} \Big). \\
\end{split}
\end{equation}
 Applying It\^o formula to $|x(t)-D(x_t)|^2$, we have:
\begin{equation}
\begin{split}
 \rvert x(t)  - D ( x_{t})\rvert^{2}  &=  \rvert x(0) - D ( \xi) \rvert^{2} +  \int_{0}^{t}  \Big [ 2[x(s) - D (x_{s})] ^{T} b( x_{s }) + \rvert \sigma (x_{s}) \rvert^{2} \Big] ds \\
& + 2\int_{0}^{s}  [x(s) - D (x_{s})] ^{T} \sigma (x_{s}) dw(s),  \\
\end{split}
\end{equation}
taking expectation on both sides,  we get
\begin{equation}\label{e}
\begin{split}
\mathbb{E}  \Big ( \sup_{ 0 < s\le t}  \rvert x(s)  - D ( x_{s})\rvert^{2} \Big)&  = \mathbb{E} \rvert x(0) - D ( \xi) \rvert^{2} \\
& +  \mathbb{E} \Big ( \sup_{ 0 <  s \le t}\int_{0}^{s} \Big [  2[x(u) - D (x_{u})] ^{T} b( x_{u }) + \rvert \sigma (x_{u}) \rvert^{2} \Big] du \Big)\\
& + 2\mathbb{E} \Big ( \sup_{ 0 <  s \le t}\int_{0}^{s} [x(u) - D (x_{u})] ^{T} \sigma (x_{u}) dw(u) \Big). \\
\end{split}
\end{equation}
Note that
\begin{equation}\label{e1}
\begin{split}
2\mathbb{E} \Big (  \sup_{ 0 <  s \le t}  \int_{0}^{s}  [x(u) - D (x_{u})] ^{T}& \sigma (x_{u}) dw(u) \Big)  \leq 8\sqrt{2} \mathbb{E} \Big ( \int_{0}^{t} \rvert [x(s) - D (x_{s})] ^{2}\rvert \sigma (x_{s})\rvert^{2} ds \Big) ^{\frac{1}{2}}\\
& \leq 12 \mathbb{E} \Big ( \int_{0}^{t}  \rvert x(s) - D (x_{s}) \rvert ^{2}  \rvert \sigma (x_{s})\rvert^{2} ds \Big)^{\frac{1}{2}}\\
& \leq 12 \mathbb{E} \Big[ \Big ( \sup_{ 0 <  s \le t}   \rvert x(s) - D (x_{s}) \rvert ^{2} \Big)^{\frac{1}{2}} \Big( \int_{0}^{t} \rvert \sigma (x_{s})\rvert^{2} ds \Big)^{\frac{1}{2}} \Big]\\
& \leq \dfrac{1}{2} \mathbb{E} \Big ( \sup_{ 0 <  s \le t}  \rvert x(s) - D (x_{s}) \rvert ^{2} \Big)+ 72   \mathbb{E}  \int_{0}^{t} \rvert \sigma (x_{s})\rvert^{2} ds.\\
\end{split}
\end{equation}
Substitute \eqref{e1} into \eqref{e} and using Lemma \ref{4.3}, yields
\begin{equation}\label{e2}
\begin{split}
\mathbb{E} & \Big ( \sup_{ 0 < s\le t}  \rvert x(s)  - D ( x_{s})\rvert^{2} \Big)  \le  (1+k)^{2} \mathbb{E}\lVert \xi \lVert^{2}_{r} +  \mathbb{E} \Big ( \sup_{ 0 <  s \le t}\int_{0}^{s} \Big [  2[x(u) - D (x_{u})] ^{T} b( x_{u })\\
& + \rvert \sigma (x_{u}) \rvert^{2} \Big] du \Big) + \dfrac{1}{2} \mathbb{E} \Big ( \sup_{ 0 <  s \le t}  \rvert x(s) - D (x_{s}) \rvert ^{2} \Big)+ 72   \mathbb{E}  \int_{0}^{t} \rvert \sigma (x_{s})\rvert^{2} ds.\\
\end{split}
\end{equation}
Applying assumption \noindent{\bf (A2)}, one have
\begin{equation}\label{e3}
\begin{split}
\mathbb{E}  \Big ( \sup_{ 0 < s\le t}&  \rvert x(s)  - D ( x_{s})\rvert^{2} \Big) \le (1+k)^{2} \mathbb{E}\lVert \xi \lVert^{2}_{r} \\
&\qquad+   \dfrac{1}{2} \mathbb{E} \Big ( \sup_{ 0 <  s \le t}  \rvert x(s) - D (x_{s}) \rvert ^{2} \Big)+ 73 L  \mathbb{E}  \int_{0}^{t} (1+\norm{x_{s}}) ds\\
& =(1+k)^{2} \mathbb{E}\lVert \xi \lVert^{2}_{r}+   \dfrac{1}{2} \mathbb{E} \Big ( \sup_{ 0 <  s \le t}  \rvert x(s) - D (x_{s}) \rvert ^{2} \Big)+73LT+73 L  \mathbb{E}  \int_{0}^{t} \norm{x_{s}} ds,
\end{split}
\end{equation}
this implies to
\begin{equation}\label{e4}
\begin{split}
\mathbb{E} & \Big ( \sup_{ 0 < s\le t}  \rvert x(s)  - D ( x_{s})\rvert^{2} \Big)  \le 2(1+k)^{2} \mathbb{E}\lVert
 \xi \lVert^{2}_{r} + 146 LT+ 146 L  \mathbb{E}  \int_{0}^{t} \norm{x_{s}} ds.
\end{split}
\end{equation}
Substituting \eqref{e4} in \eqref{ee}, we obtain,
\begin{equation}\label{e9}
\begin{split}
\mathbb{E}\big(\sup_{0\leq s\leq
	t} \rvert x(s) \rvert^{2}\big)&\leq  2k^{2}E\|x_{t}\|_r^2+ 4(1+k)^{2} \mathbb{E}\lVert
\xi \lVert^{2}_{r} + 292 LT+ 292 L  \mathbb{E}  \int_{0}^{t} \norm{x_{s}} ds. \\
\end{split}
\end{equation}
Now, correspond to the definition of the norm $ \rVert \cdot \rVert _{r} $, we have:
\begin{equation}\label{5.8}
\begin{split}
\mathbb{E} \rVert x_{t} \rVert _{r}^{2} & = \mathbb{E} \Big( \sup_{ - \infty < \theta \le 0} e^{r \theta} \rvert x_{t} (\theta) \rvert \Big) ^{2} \\
&\leq  e ^{- 2r t} \mathbb{E} \rVert  \xi \rVert^{2}_{r} +  \mathbb{E} \Big( \sup_{ 0 < s \le t} \rvert x ( s ) \rvert^{2} \Big).
\end{split}
\end{equation}
by substituting this into \eqref{e9}, we get
\begin{equation}
\begin{split}
\mathbb{E}\big(\sup_{0\leq s\leq
	t} \rvert x(s) \rvert^{2}\big)&\leq  2k^{2} \mathbb{E} \rVert  \xi \rVert^{2}_{r} + 4(1+k)^{2} \mathbb{E}\lVert\xi \lVert^{2}_{r} + 292 LT\\
& + 292 L  \mathbb{E}  \int_{0}^{t}\Big[ e ^{-2rs} \mathbb{E} \rVert  \xi \rVert^{2}_{r} +  \mathbb{E} \Big(  \rvert x ( s ) \rvert^{2} \Big) \Big]ds, \\
\end{split}
\end{equation}
this yields to,
\begin{equation}
\begin{split}
\mathbb{E}\big( \rvert x(t) \rvert^{2}\big)&\leq \Big[\Big( 2k^{2}  + 4(1+k)^{2}+\dfrac{292L}{2r}\Big) \mathbb{E} \rVert  \xi \rVert^{2}_{r}+ 292 LT\Big] + 292 L  \mathbb{E}  \int_{0}^{t}  \mathbb{E} \Big(  \rvert x ( s ) \rvert^{2} \Big) ds, \\
\end{split}
\end{equation}
so, by the Gronwall iniquity we get,
\begin{equation}
\begin{split}
\mathbb{E}\big(\sup_{0\leq s\leq
	t} \rvert x(s) \rvert^{2}\big)&\leq C_{4}e^{C_{5}t}, \\
\end{split}
\end{equation}
where, $ C_{4}= \Big[\Big( 2k^{2}  + 4(1+k)^{2}+\dfrac{292L}{2r}\Big) \mathbb{E} \rVert  \xi \rVert^{2}_{r}+ 292 LT\Big] $ and $ C_{5}=292 L $.
The proof is complete.\hfill $\Box$\\
\begin{thm} Let assumption \noindent{\bf (A2)} and \noindent{\bf (A3)} hold. Then
		\begin{equation}\label{e8}
		\lim_{t\rightarrow\infty}\sup \frac{1}{t}log|x(t)|\leq \hat{C},
		\end{equation}
		where $\hat{C}=\dfrac{292L}{2}$.
	\end{thm}
\textbf{Proof:} From the Theorem \ref{est}, for each $m=1,2,3,...$ we get
	\begin{equation*}\mathbb{E}[\sup_{m-1\leq t\leq
		m}x(s)]\leq C_4e^{292 L m},
	\end{equation*}
	where $ C_{4}= \Big[\Big( 2k^{2}  + 4(1+k)^{2}+\dfrac{292L}{2r}\Big) \mathbb{E} \rVert  \xi \rVert^{2}_{r}+ 292 LT\Big] $. By Chebyshev inequality we further obtain that
	\begin{equation}\label{n2}
	P\{w:\sup_{m-1\leq t\leq m}|x(s)|^2>e^{(292 L +\epsilon)m}\}\leq
	C_4e^{-\epsilon m}.
	\end{equation}
	As the series $\sum_{m=1}^{\infty}C_4e^{-\epsilon m}$ is convergent,
	for almost all $w\in\Omega$, the Borel-Cantelli lemma yields that
	there exists a random integer $m_0=m_0(w)$ such that
	\begin{equation*}
	\sup_{m-1\leq t\leq m}|x(t)|^2\leq e^{(292 L +\epsilon)m},\,\,
	\textit{whenever}\,\,m\geq m_0,\end{equation*} That is, for $m-1\leq
	t\leq m$ and $m\geq m_0$, we derive
	\begin{equation*}
	|x(t)|\leq e^{\frac{1}{2}(292 L +\epsilon)m}.\end{equation*}
	Hence for almost all $w\in\Omega$ if $m-1\leq
	t\leq m$ and $m\geq m_0$, then
	\begin{equation*}
	\lim_{t\rightarrow\infty}sup\frac{1}{t}Log|x(t)|\leq \frac{1}{2}(292 L +\epsilon)m,\end{equation*} the required claim \eqref{e8}
	follows because $\epsilon>0$ is arbitrary.\hfill $\Box$\\

\section*{Acknowledgement}
I would like to thank Dr Jianhai Bao for his encouragement and kindly advice throughout this work, as well as Professor Chenggui Yuan for his supervision and remarks. This research was supported by Kufa University and the Iraqi Ministry of Higher Education and Scientific Research.

\end{document}